\documentclass[11pt, oneside]{article}   	
\usepackage{geometry}                		
\usepackage{graphicx}				
\usepackage{amssymb}
\usepackage{amsmath}
\usepackage{amsthm}

\newtheorem{conj}{Conjecture}

\title{On  a conjecture about an analogue of Tokuyama's theorem for $G_2$}
\author{Mario DeFranco}

\begin{document}
\maketitle

\begin{abstract}
We prove the conjecture of \cite{FGG} about relating sums over Littelmann patterns to the  the root system of type $G_2$, which is an analogue of Tokuyama's theorem \cite{T} for root systems of type $A_r$. We use elementary means to show that the conjecture is implied by a finite set of polynomial identities.
\end{abstract}

For each $r\geq 1$, Tokuyama's theorem \cite{T} relates a sum over Gel'fand-Cetlin patterns to the product of the Weyl character formula for $A_r$ and a polynomial we denote by $D(\textbf{x}; A_r)$. For general roots systems $R$, there is a set of inequalities called \textit{Littelmann patterns} which generalize the Gel'fand-Cetlin patterns for type $A_r$. In \cite{FGG}, the authors study a sum over Littelmann patterns for the root system of type $G_2$ and conjecture that it is equal to the product of Weyl character formula for $G_2$ and a polynomial $D(\textbf{x}; G_2)$ defined below. This formula, expressed in Conjecture \ref{conjecture}, is thus an analogue of Tokuyama's theorem for $G_2$. 

Tokuyama's theorem is an important part in studying Weyl Group Multiple Dirichlet series, which are series defined using root systems. These series are also related to $p$-adic Whittaker functions. This connection is described in \cite{B}.  

We prove Conjecture \ref{conjecture} by expressing both sides as polynomials in four indeterminates whose coefficients are rational functions. We then show that the coefficients are equal. 

We define the terms necessary to state Conjecture \ref{conjecture}. The set $\{ \alpha_1,\alpha_2\}$ is a choice of simple roots for $G_2$, with $\alpha_2$ being the longer root. Let $W$ be the Weyl group for $G_2$, $\Lambda_W$ the weight lattice, and $\mathbb{C}[\Lambda_W]$ the associated ring of Laurent polynomials. Let $\varpi_1,\varpi_2 \in \Lambda_W$ be the fundamental weights for $G_2$, and $\rho$ be the half-sum of positive roots. For a dominant weight $\theta \in \Lambda_W$, the weight $\theta + \rho = \ell_1 \varpi_1 + \ell_2 \varpi_2$,  for some positive integers $\ell_1$ and $\ell_2$. The variables $x$ and $y$ are the indeterminates of $\mathbb{C}[\Lambda_W]$. Let $\textbf{x}^{m_1\alpha_1+m_2\alpha_2} = x^{m_1} y^{m_2}$. $D(\textbf{x})$ denotes 
\[
D(\textbf{x})=\prod_{\alpha>0} (1-\frac{\textbf{x}^\alpha}{q})=(1-\frac{x}{q})(1-\frac{y}{q})(1-\frac{x y}{q})(1-\frac{x^2 y}{q})(1-\frac{x^3 y}{q})(1-\frac{x^3 y^2}{q}).
\]

 The set $\mathcal{B}(\theta+\rho)$ is the set of Littlemann patterns, which are 6-tuples $\pi=(a,b,c,d,e,f)$, where $a,b,c,d,e,f$ are non-negative integers that satisfy the following Littelmann inequalities. 
\begin{enumerate} 
\item $0 \leq \, f \leq \ell_2 +a - 2b+ c-2d+e$  \label{f} 
\item  $b \leq \, a \leq \ell_1 +3b - 2c+ 3d-2e$  \label{a} 
\item  $\frac{c}{2} \leq \, b \leq \ell_2 +c-2d+e$ \label{b}  
\item $ 2d \leq \, c \leq \ell_1 +3d-2e$ \label{c} 
\item  $e \leq \, d \leq \ell_2 +e$ \label{d}
\item   $0 \leq \, e \leq \ell_1$ \label{littelmann inequalities} 
\end{enumerate}

  It is straightforward to check that inequalities $i+1$ through $6$ imply that the lower bound of inequality $i$ is less than or equal to the upper bound of inequality $i$. Following the terminology of \cite{FGG}, we say that an entry $u$ of $\pi$ is ``circled", denoted by $u^\circ$, if $u$ attains its lower bound; e.g., $f$ is circled if $f=0$, $a$ is circled if $a=b$, $b$ is circled if $b=\frac{c}{2}$, etc. We say that $u$ is ``boxed", denoted by $\underline{u}$, if $u$ attains its upper bound in its inequality; e.g., $f$ is boxed if $f = \ell_2 +a - 2b+ c-2d+e$, $a$ is boxed if $a=  \ell_1 +3b - 2c+ 3d-2e$, $b$ is boxed if $b= \ell_2 +c-2d+e$, etc. Then for a Littlemann pattern $\pi$, \cite{FGG} define $\hat{H}(\pi) \in \mathbb{Z}[q^{-1}]$. Their conjecture is then 

\begin{conj} \label{conjecture}
\[
\sum_{\pi \in \mathcal{B}(\theta+\rho)}   \hat{H}(\pi) \,  x^{a+c+e} y^{b+d+f} = \textbf{x}^{-w_\ell(\theta+\rho)} D(\textbf{x})  \frac{\sum_{w \in W} \mathrm{sgn}(w) \textbf{x}^{w(\theta+\rho)}}{\prod_{\alpha>0} (1-x^{\alpha})}
\]
\end{conj}
 We note that to define $\hat{H}(\pi)$, \cite{FGG} use a definition depending on whether $\pi$ is generic or one of twenty special cases.  
We will express the coefficient $\hat{H}(\pi)$ as 
\begin{equation} 
 \hat{H}(\pi)=H_{\mathrm{std}}(\pi) + H_{\mathrm{adj}}(\pi) 
 \end{equation}
where the the generic part is encompassed by the ``standard" term $H_{\mathrm{std}}(\pi)$ and the special cases by the ``adjusted" term $H_{\mathrm{adj}}(\pi)$. This allows us to consolidate the special cases, to simplify their characterization, and also to simplify the values of $H_{\mathrm{adj}}(\pi)$. We also see that the after expressing 
\[
\sum_{\pi \in \mathcal{B}(\theta+\rho)}   H_{\text{adj}}(\pi) \,  x^{a+c+e} y^{b+d+f} 
\]   
as a polynomial, the coefficients have a factored form similar to that of $D(\textbf{x})$.

Now we define $H_\text{std}(\pi)$ and $H_\text{adj}(\pi)$. The coefficient $H_{\mathrm{std}}(\pi)$ denotes the ``standard contribution'' for the 6-tuple $\pi = (a,b,c,d,e,f)$ defined by 
  \[
 H_\mathrm{std}(\pi) = h(a)h(b)h(c)h(d)h(e)h(f)  
  \]
  where 
  \[
  h(u) = \begin{cases} 1-\frac{1}{q}, u \, \text{is neither boxed nor circled}  \\ 
   				\frac{-1}{q},  u \, \text{is boxed} \\ 
				1,  u \, \text{is circled} \\ 
				0, u \, \text{is both boxed and circled}.\\  						
  \end{cases}
  \]

The ``adjusted contribution" $H_{\text{adj}}(\pi)$ is defined in general to be 0 unless $\pi=(a,b,c,d,e,f)$ satisfies certain conditions. The first condition is that  
$\pi$ has what \cite{FGG} define to be a ``bad middle", which means $b=d+1$ and $c= 2d+1$. Therefore the Littelmann inequalities for $\pi$ with bad middles become 
\begin{enumerate}
\item $0 \leq \, f \leq \ell_2 +a -2d+e-1$  \label{f} 
\item  $d+1 \leq \, a \leq \ell_1+ 2d-2e+1$  \label{a} 
\item  $d+.5 \leq \, d+1 \leq \ell_2 +e+1$ \label{b}  
\item $ 2d \leq \, 2d+1 \leq \ell_1 +3d-2e$ \label{c} 
\item  $e \leq \, d \leq \ell_2 +e$ \label{d}
\item   $0 \leq \, e \leq \ell_1$\label{e}
\end{enumerate}
Thus such $\pi$ are determined by the values of $e,d,a$ and $f$. The definitions for circling and boxing the entries of $\pi$ still hold. 


 

We define for any $\pi \in \mathcal{B}(\theta+\rho)$
\[
H_\text{adj}(\pi) = \hat{H}(\pi) -H_\mathrm{std}(\pi)
\] 
where $\hat{H}(\pi)$ is defined by \cite{FGG} according to some twenty cases.  Let $\pi'=(a,b,c,d,e)$, and we set 
\[
H_\text{adj}(\pi) = H_\text{adj}(\pi') h(f).
\]  
By calculating $\hat{H}(\pi)$ and $H_\mathrm{std}(\pi)$ in each of these cases, we can determine $H_\text{adj}(\pi)$. We see that the values of $H_\text{adj}(\pi')$ become more concise than those for $\hat{H}(\pi')$ given in \cite{FGG} and that the twenty cases are consolidated to the following definition. 
\[
H_\text{adj}(\pi') = \begin{cases}
			 \frac{(1-1/q)}{q}, (e^\circ,d^\circ,a^\circ) \\ 
			 -\frac{(1-1/q)}{q^2}, (e^\circ, d \text{ or } d^\circ,\underline{a}) \\ 
			 \frac{(1-1/q)}{q^3}, (e \text{ or } e^\circ, \underline{d}, a = 2d+1-e)\\
			  \frac{(1-1/q)^2}{q},  (e,d^\circ,a^\circ) ,  (e^\circ,d^\circ,a),  (e^\circ,d,a^\circ),  (e,d^\circ,a),  (e,d,a^\circ) \\   
			  -\frac{(1-1/q)^2}{q^2},(e \text{ or } e^\circ, \underline{d}, a) \text{ such that } a \neq 2d+1-e \\ 
			  \frac{(1-1/q)^3}{q}, (e,d,a) \\  
			  \frac{(1-1/q)^3}{q}, (e^\circ, d, a) \text{ such that } a \neq 2d+1-e \\ 
			  \frac{(1-1/q)}{q}((1-1/q)^2 -1/q ),    (e^\circ, d, a = 2d+1-e)	.		  
			   \end{cases}
\]
This means, for example, that if $\pi' = (1,1, 1, 0,0)$, then that means $\pi'$ has a bad middle with $e,d$ and $a$ circled (because we assume $\ell_1, \ell_2>0$), so 
\[
H_\text{adj}((1,1,1,0,0)) =   \frac{(1-1/q)}{q}.
\]
We see that the definition of $H_{\text{adj}}(\pi)$ depends only on the circling and boxing of $e,d$ and $a$ and whether $a = 2d+1-e$. 

Now we can prove conjecture \eqref{conjecture}.

\begin{proof}
The strategy of the proof is to express  
\begin{equation} \label{littelmann sum}
\sum_{\pi \in \mathcal{B}(\theta+\rho)}   \hat{H}(\pi) \,  x^{a+c+e} y^{b+d+f} 
\end{equation}
as a rational function in $x,y$ and $q^{-1}$. This rational function depends on the numbers $\ell_1$ and $\ell_2$, which only appear as exponents of $x$ and $y$ in the numerator of the rational function. We therefore interpret this rational function as a polynomial, say $P_H$, in the four indeterminates 
\begin{equation} \label{four indeterminates}
x^{\ell_1}, y^{\ell_1}, x^{\ell_2}, y^{\ell_2} 
\end{equation}
whose coefficients we prove will be of the form 
\begin{equation} \label{coefficient form}
\frac{p_1(x,y,q^{-1})}{p_2(x,y)}
\end{equation} 
where $p_1$ and $p_2$ are polynomials. Now the right side of \eqref{conjecture} is also a polynomial, say $P_W$, in the four indeterminates \eqref{four indeterminates} with coefficients of the form \eqref{coefficient form}. Therefore equality of \eqref{conjecture} can be established by equating the coefficients of the two polynomials $P_H$ and $P_W$. The polynomial $P_W$ has 12 terms, as there are 12 elements in the Weyl group $W$ and the coefficients are of the form 
\[
\text{sgn}(w) T(\textbf{x})
\]. 
 We denote the multi-degree of the term 
\[
(x^{\ell_1})^{m_1} (y^{\ell_1})^{n_1}(x^{\ell_2})^{m_2}  (y^{\ell_2})^{n_2}
\]
 by
 \[
( (m_1,n_1), (m_2,n_2)).
 \]
Define 
\[
T(\textbf{x}) = \frac{(1-q^{-1}x)(1-q^{-1}y)(1-q^{-1}xy)(1-q^{-1}x^2y)(1-q^{-1}x^3y)(1-q^{-1}x^3y^2)}{(1-x)(1-y)(1-xy)(1-x^2y)(1-x^3y)(1-x^3y^2)}.
\]
Then the twelve multi-degrees of $P_W$ and the coefficients are given in Table \ref{tab:tableW}.
\begin{table}[h!]
  \begin{center}
    \caption{Terms for $P_W$}
    \label{tab:tableW}
    \begin{tabular}{l|c|r} 
      \textbf{multi-degree} & \textbf{coefficient}\\
      \hline
  ((1,0),(0,0)) &$-T(\textbf{x})$\\
 ((1,1),(0,1)) & $T(\textbf{x})$\\
  ((0,0),(0,0))  & $T(\textbf{x})$\\
  ((0,0),(0,1)) & $-T(\textbf{x})$ \\ 
  ((1,0),(3,1))  & $-T(\textbf{x})$\\
  ((3,1),(3,1))  & $T(\textbf{x})$ \\
  ((3,1),(6,3))& $-T(\textbf{x})$ \\
  ((4,2),(6,3)) & $T(\textbf{x})$\\
  ((1,1),(3,3)) & $-T(\textbf{x})$ \\
  ((3,2),(3,3)) & $-T(\textbf{x}) $\\
  ((3,2),(6,4)) & $T(\textbf{x})$\\
  ((4,2),(6,4)) & $-T(\textbf{x})$ \\.

 	\end{tabular}
  \end{center}
\end{table}

We show how to express \eqref{littelmann sum} as a polynomial in the indeterminates \eqref{four indeterminates}. As $\hat{H}(\pi) = H_\text{std}(\pi)+ H_\text{adj}(\pi)$, 
we compute separately the two sums 
\begin{equation} \label{littelmann sum H std}
\sum_{\pi \in \mathcal{B}(\theta+\rho)}   H_\text{std}(\pi) \,  x^{a+c+e} y^{b+d+f} 
\end{equation}
and 
\begin{equation} \label{littelmann sum H adj}
\sum_{\pi \in \mathcal{B}(\theta+\rho)}   H_\text{adj}(\pi) \,  x^{a+c+e} y^{b+d+f}. 
\end{equation}
We first compute \eqref{littelmann sum H std}. We sum over the six indices in the order $f, a, b, c, d, e$. We write \eqref{littelmann sum H std} as 
\begin{equation} \label{nested abcdef sum}
 \sum_e h(e)x^e  \left( \sum_d h(d)y^d  \left( \sum_c h(c)x^c \left( \sum_b h(b)y^b  \left( \sum_a h(a)x^a \left( \sum_f h(f)y^f \right) \right) \right) \right) \right)
\end{equation}
where the indices are over the Littelmann inequalities. Thus $f$ ranges from 0 to $\ell_2 +a - 2b+ c-2d+e$, $a$ ranges from $b$ to $\ell_1 +3b - 2c+ 3d-2e$, $b$ ranges from $\lceil c/2 \rceil$ to $\ell_2 +c-2d+e$, etc. We evaluate these sums in the following way. 

Let $u$ be an entry of $\pi$ and $L$, $U$ the lower and upper bounds in the Littelmann inequality for $u$. Then, if $U >L$, 
\begin{align} \label{u sum}
\sum_{L \leq u \leq U} h(u)X^u  &= X^L + (1-q^{-1})\frac{X^{L+1} -X^{U}}{1-X} -q^{-1} X^{U} \\ 
								           &=  \frac{(1-q^{-1} X)}{1-X}(X^L-X^U).
\end{align}
This equation also holds when $U=L$, as both sides are 0. Thus \eqref{u sum} is a polynomial in the indeterminates \eqref{four indeterminates} with coefficients of the form \eqref{coefficient form}.   
The only issue in evaluating these sums is that the lower bound for $b$ is $\lceil c/2 \rceil$. To evaluate the sum over $b$ we make use of characteristic functions $\textbf{1}_0$, where 
\[
\textbf{1}_0(n) = \begin{cases} 
			1, n \equiv 0 \mod 2\\ 
			0, n \equiv 1 \mod 2. \\			
			\end{cases}
\]
We then have 
\begin{align}
\sum_{\lceil c/2 \rceil \leq b\leq U } h(b) X^b = (X^{\frac{c}{2}}+ \frac{X^{\frac{c}{2}+1}(1-q^{-1})}{1-X})\textbf{1}_0(c) \nonumber \\ 
								 + (\frac{X^{\frac{c+1}{2}}(1-q^{-1})}{1-X})\textbf{1}_0(c+1) \nonumber \\ 
								- \frac{X^U(1-q^{-1})}{1-X} - q^{-1}X^U. \label{b sum}
\end{align}
This equation \eqref{b sum} also holds when $U=\lceil c/2 \rceil$. Therefore \eqref{b sum} leads us to evaluate sums of the form  
\[
\sum_{L\leq u \leq U} h(u)Y^u  X^{\frac{C_1+C_2u}{2}} \textbf{1}_0 (C_1+C_2u) 
\]
where $C_1$ and $C_2$ are odd integers and $U\geq L$. We obtain  
\begin{align} \label{u sum characteristic} 
\sum_{L \leq u \leq U} h(u)Y^u X^{\frac{C_1+C_2u}{2}} \textbf{1}_0 (C_1+C_2u)  &=  (Y^L X^{\frac{C_1+C_2L}{2}} +\frac{(1-q^{-1})Y^{L+2}X^{\frac{C_1+C_2(L+2)}{2}} }{1-Y^2X^{C_2}})\textbf{1}_0 (C_1+C_2L) \\
								           & + \frac{(1-q^{-1})Y^{L+1}X^{\frac{C_1+C_2(L+1)}{2}}}{1-Y^2X^{C_2}}\textbf{1}_0 (C_1+C_2(L+1))\\ 
								           &-(\frac{(1-q^{-1})Y^{U}X^{\frac{C_1+C_2U}{2}}}{1-Y^2X^{C_2}} + q^{-1}Y^UX^{\frac{C_1+C_2U}{2}})\textbf{1}_0 (C_1+C_2 U) \\ 
								           & - \frac{(1-q^{-1})Y^{U+1}X^{\frac{C_1+C_2(U+1)}{2}}}{1-Y^2X^{C_2}}\textbf{1}_0 (C_1+C_2(U+1)).
\end{align}
As the sums over $c,d$ and $e$ all have integral upper and lower bounds in their Littelmann inequalities, equation \eqref{u sum characteristic} suffices to evaluate these sums.
In this way we can express 
\begin{equation} \label{Hstd sum}
\sum_{\pi \in \mathcal{B}(\theta+\rho)} H_{\mathrm{std}}(\pi) \,  x^{a+c+e} y^{b+d+f} 
\end{equation}
as a finite sum of terms of the form 
\begin{equation} \label{characteristic term form}
\frac{P_1(x,y,q^{-1})}{P_2(x,y)}(x^{\ell_1})^{n_1} (x^{\ell_2})^{n_2} (y^{\ell_1})^{n_3} (y^{\ell_2})^{n_4} \textbf{1}_0(A_1\ell_1 + A_2 \ell_2 + A_3)
\end{equation}
where $n_i$ are non-negative integers and  $A_i$ are integers. 

Likewise we can express 
\begin{equation} \label{Hadj sum}
\sum_{\pi \in \mathcal{B}(\theta+\rho)}  H_{\mathrm{adj}}(\pi) \,  x^{a+c+e} y^{b+d+f}
\end{equation}
as a finite sum of terms of the form 
\begin{equation} \label{term form}
\frac{P_1(x,y,q^{-1})}{P_2(x,y)}(x^{\ell_1})^{n_1} (x^{\ell_2})^{n_2} (y^{\ell_1})^{n_3} (y^{\ell_2})^{n_4}.
\end{equation}
Because $c$ is always odd in the cases for $H_\text{adj}$, we do not need the characteristic functions.

To equate the coefficients with the right side of \eqref{conjecture}, we have to specify parities for $\ell_1$ and $\ell_2$ to render \eqref{Hstd sum} a true polynomial without characteristic functions. That is, we set   
\[
\ell_1 = 2m_1+\epsilon_1, \, \ell_1 = 2m_2+\epsilon_2
\]
for a choice of $\epsilon_1, \epsilon_2 \in \{0,1 \}$. Then \eqref{Hstd sum} becomes a polynomial in the indeterminates
\[
x^{2n_1}, x^{2n_2}, y^{2n_1}, y^{2n_2} 
\] 
with coefficients of the form \eqref{coefficient form}.

We calculate that the sum of the standard terms \eqref{Hstd sum} is equal to a sum of 544 terms of the form \eqref{characteristic term form} and that the sum for adjusted terms \eqref{Hadj sum} is a sum of 106 terms of the form \eqref{term form}. We then consider the set of multi-degrees that occurs in these sums. We denote the multi-degree of the term 
\[
\frac{P_1(x,y,q^{-1})}{P_2(x,y)}(x^{\ell_1})^{m_1} (y^{\ell_1})^{n_1}(x^{\ell_2})^{m_2}  (y^{\ell_2})^{n_2} \textbf{1}_0(A_1\ell_1 + A_2 \ell_2 + A_3)
\]
 by
 \[
( (m_1,n_1), (m_2,n_2)).
 \]
 There are 33 distinct multi-degrees that occur from the standard terms \eqref{Hstd sum}, and 14 distrinct multi-degrees that come from the adjusted terms \eqref{Hadj sum}. The union of these sets contains 35 distinct multi-degrees. When we combine like terms for the standard terms, there are 18 multi-degrees with non-zero coefficients, and when we combine like terms for the adjusted terms, there are 10 multi-degrees with non-zero coefficients. We present the multi-degrees and coefficients with $\epsilon_1 = \epsilon_2=0$ for the standard terms in Table \ref{tab:tablestd} and for the adjusted terms in Table \ref{tab:tableadj}. To express these coefficients, we define 
 \begin{align} 
 T_1(\textbf{x}) &= \frac{(1-q^{-1})(1-q^{-1}x)(1-q^{-1}y)(1-q^{-1}x^3y^2)}{(1-x)(1-y)(1-x^4y^2)(1-x^3y^2)}\\ 
 T_2(\textbf{x}) &= \frac{(1-q^{-1})(1-q^{-1}y)(1-q^{-1}xy)(1-q^{-1}x^3y)}{(1-y)(1-xy)(1-x^4y^2)(1-x^3y)} \\ 
 T_3(\textbf{x}) &= \frac{(1-q^{-1})^2(1-q^{-1}y)(1-q^{-1}x^4y^2)}{(1-x)(1-xy)(1+x^2y)(1-x^3y)(1-x^3y^2)}
 \end{align}
 in addition to $ T(\textbf{x}) $ defined above.
\vspace{1 in}

\begin{table}[h!]
  \begin{center}
    \caption{Standard terms with  $\epsilon_1 = \epsilon_2=0$}
    \label{tab:tablestd}
    \begin{tabular}{l|c|r} 
      \textbf{multi-degree} & \textbf{coefficient}\\
      \hline
  ((1,0),(0,0)) &$-T(\textbf{x})+q^{-1}x^2y T_1(\textbf{x})$\\
 ((1,1),(0,1)) & $T(\textbf{x})-q^{-1}x^2y  T_2(\textbf{x})$\\
  ((1,0),(4,2)) & $-q^{-1}x^2y T_1(\textbf{x}) $\\
   ((1,1),(4,3)) & $q^{-1}x^2y T_2(\textbf{x})$ \\
  ((0,0),(3,2))  &  $q^{-1}x^3y T_3(\textbf{x}) $\\
  ((4,2),(3,2))  &  $-q^{-1}x^3y T_3(\textbf{x}) $\\ 
  ((4,2),(4,2))  & $q^{-1}x^3y T_1(\textbf{x}) $\\
  ((4,2),(4,3))  & $-q^{-1}x^3y T_2(\textbf{x}) $\\
  ((0,0),(0,0))  & $T(\textbf{x})-q^{-1}x^2y  T_1(\textbf{x})$\\
  ((0,0),(0,1)) & $-T(\textbf{x})+q^{-1}x^2y  T_2(\textbf{x})$ \\ 
  ((1,0),(3,1))  & $-T(\textbf{x})$\\
  ((3,1),(3,1))  & $T(\textbf{x})$ \\
  ((3,1),(6,3))& $-T(\textbf{x})$ \\
  ((4,2),(6,3)) & $T(\textbf{x})$\\
  ((1,1),(3,3)) & $-T(\textbf{x})$ \\
  ((3,2),(3,3)) & $-T(\textbf{x}) $\\
  ((3,2),(6,4)) & $T(\textbf{x})$\\
  ((4,2),(6,4)) & $-T(\textbf{x})$ \\

 	\end{tabular}
  \end{center}
\end{table}

\begin{table}[h!]
  \begin{center}
    \caption{Adjusted terms}
    \label{tab:tableadj}
    \begin{tabular}{l|c|r} 
      \textbf{multi-degree} & \textbf{coefficient}\\
      \hline
((1,0),(0,0)) & $-q^{-1}x^2y T_1(\textbf{x})$ \\
 ((1,1),(0,1)) &$q^{-1}x^2y  T_2(\textbf{x})$ \\
((1,0),(4,2))  &$q^{-1}x^2y  T_1(\textbf{x})$\\
((1,1),(4,3))  &$-q^{-1}x^3y  T_2(\textbf{x})$\\
((0,0),(3,2))  &$-q^{-1}x^3y  T_3(\textbf{x})$\\
((4,2),(3,2)   &$q^{-1}x^2y  T_3(\textbf{x})$ \\
 ((4,2),(4,2))  &$-q^{-1}x^2y  T_1(\textbf{x})$\\
 ((4,2),(4,3)) &$q^{-1}x^2y  T_2(\textbf{x})$\\
((0,0),(0,0)) & $q^{-1}x^2y  T_1(\textbf{x})$\\
((0,0),(0,1))  & $-q^{-1}x^2y  T_2(\textbf{x})$\\
 	\end{tabular}
  \end{center}
\end{table}
 
Now, sums of the form \eqref{u sum} do evaluate to polynomials in $X^L$ and $X^U$, but sums of the form \eqref{b sum}  and \eqref{u sum} in general do not evaluate to such polynomials. For example, there is no polynomial $P(Z)$ such that
\[
\sum_{\lceil c/2 \rceil \leq b\leq U } h(b) X^b = P(X^{c/2}).
\]
However, the entire sum 
\[
\sum_{\pi \in \mathcal{B}(\theta+\rho)}   H_\text{std}(\pi) \,  x^{a+c+e} y^{b+d+f} 
\]
 is equal to a polynomial in the indeterminates \eqref{four indeterminates}. We verify this by computing the Table \ref{tab:tablestd} for all four combinations of $\epsilon_1$ and $\epsilon_2$ and seeing that the coefficients agree.  
 
Now we add the coefficients for the standard terms and adjusted terms for each multi-degree and see that they add up to 
\[
\text{sgn}(w) T(\textbf{x})
\]
 which gives us the right side of \eqref{conjecture} and proves the result. 
 \end{proof}


\begin{thebibliography}{9}


\bibitem{B} Bump, D.: Introduction: multiple Dirichlet series, Multiple Dirichlet series, L-functions and automorphic forms, Progr. Math., vol. 300, pp. 1?36. Birkhäuser/Springer, New York (2012)

  \bibitem{FGG}
 Friedlander, H., Gaudet, L. \& Gunnells, P.E.: Crystal graphs, Tokuyama's theorem, and the Gindikin-Karpelevic formula for $G_2$. J Algebr Comb (2015) 41: 1089. https://doi.org/10.1007/s10801-014-0567-9 
 
 \bibitem{L}
Littelmann, P.: Cones, crystals, and patterns. Transform. Groups 3(2), 145-179 (1998)
 
 \bibitem{T}
Tokuyama, T.: A generating function of strict Gelfand patterns and some formulas on characters of general linear groups. J. Math. Soc. Jpn. 40(4), 671-685 (1988)
 \end{thebibliography}
\end{document}